\documentclass[12pt,a4paper,leqno]{amsart}

\usepackage{anysize}
\usepackage[utf8]{inputenc}
\usepackage{amsmath}
\usepackage{amssymb}
\usepackage{array}
\usepackage[unicode]{hyperref}

\marginsize{3cm}{3cm}{3.5cm}{3.5cm}

\newtheorem{theorem}{Theorem}

\begin{document}

\title[The fastest series for $1/\pi$ due to Ramanujan]{The fastest series for $1/\pi$ due to Ramanujan. \\ Proofs from modular polynomials}
\author{Jesús Guillera} 
\address{University of Zaragoza, Department of mathematics, 50009 Zaragoza (Spain)}
\email{}

\begin{abstract}
First we give general formulas for proving real or complex Ramanujan series for $1/\pi$. Then, as an example, we apply them for providing complete proofs of the fastest series for $1/\pi$ due to Ramanujan using Russell and Weber modular polynomials. We recommend the reader to use the Maple code in \cite{Gui-maple-program} for automatically proving any Ramanujan-type series for $1/\pi$.
\end{abstract}

\maketitle

\section{Introduction}
First we give formulas for proving real or complex Ramanujan series for $1/\pi$. Then, as an example, we apply them for obtaining complete proofs, using Russell \cite{Russell} and Weber modular polynomials \cite{Weber-polynomial}, of the fastest formula (8 digits per term) due to Ramanujan \cite[eq. 44]{Ramanujan} 
\begin{equation}\label{Ramanujan-1103}
\sum_{n=0}^{\infty} \frac{\left(\frac12\right)_n\left(\frac14\right)_n \left(\frac34 \right)_n}{(1)_n^3} (26390n+1103) \left( \frac{1}{99^4} \right)^n = \frac{9801\sqrt 2}{4 \pi}. 
\end{equation}
J. and P. Borwein were the first in trying to understand how Ramanujan arrived at his formula \cite[Chapter 5]{Borwein-AGM}. They could not complete the proof because the difficult step uses the value of a Ramanujan's multiplier that they only could evaluate numerically. In addition, they commented that in their opinion Ramanujan probably did the same. Complete and rigorous proofs are \cite{Zhao-RAMA-modular} and \cite[Sect. 6]{Wronskians-Theta}. Some interesting surveys on Ramanujan series are \cite{BaBeCh}, \cite{Rama-Zudilin}, \cite{Rama-Sato}, \cite{Cohen-Gui} and \cite{Gui-Rama-Pi-survey}. 
\par The theory that we are going to develop is only valid for hypergeometric series (the types considered by Ramanujan), that is, for the levels $1, 2, 3, 4$, and is completely rigorous. Certainly we make some numerical evaluations to know the degree of the corresponding modular equation, but then from the modular equation we obtain the exact value of the degree confirming the hypothesis. In addition we believe that our method could be extended to higher levels. Recall that In 2002 Takeshi Sato surprised the mathematical community by showing in a conference a non-hypergeometric series for $1/\pi$ (see the Wikipedia entry on Ramanujan--Sato series).
\par We could not omit an explanation of the Maple program in \cite{Gui-maple-program} which proves automatically any Ramanujan-type series for $1/\pi$. This is done in the last section. In which follows we develop our theory.

\subsection{Legendre's relation}

We will use the following version: Let
\[
\ell=4 \sin^2 \frac{\pi}{s}, 
\]
where $s=6,\, 4, \, 3, \, 2$, or equivalently $\ell=1, \, 2, \,3 , \, 4$, and
\[
F_\ell(x)={}_2F_1\biggl(\begin{matrix} \frac1s, \, 1-\frac1s \\ 1 \end{matrix}\biggm| x \biggr), \quad G_{\ell}(x) = x \frac{d F_{\ell}(x)}{dx}, 
\]
then we have
\[
\frac{2 \alpha}{\sqrt{\ell}} F_{\ell}(\alpha) G_{\ell}(\beta) + 
\frac{2 \beta}{\sqrt{\ell}} F_{\ell}(\beta) G_{\ell}(\alpha) = \frac{1}{\pi}.
\]
when $\beta=1-\alpha$, which explains why $\pi$ appears in the Ramanujan-type series.

\subsection{Clausen's relation}

We will use the following version: We have
\[
{}_2F_1\biggl(\begin{matrix} \frac1s, \, 1-\frac1s \\ 1 \end{matrix}\biggm| x \biggr)^2 =
{}_3F_2\biggl(\begin{matrix} \frac12, \, \frac1s, \, 1-\frac1s \\ 1, \, 1 \end{matrix}\biggm| 4 x(1-x) \biggr)
\]
for ${\rm Re}(x)<1/2$.

\subsection{The nome $q$ and the mirror map}

If we write
\[
F_\ell^2(x) =\frac{q}{x(1-x)} \frac{d x}{dq},
\]
we get
\[
q= B_{\ell} \cdot \exp \int \frac{dx}{F_{\ell}(x)^2 x (1-x)},
\]
where $B_{\ell}$ is the minimal constant such that the expansion of $q$ in powers of $x$ has integer coefficients. The inverse function $x_{\ell}(q)$ is the mirror map. If we let  $\beta=x_{\ell}(q)$ and $\alpha=x_{\ell}(q^d)$, then 
\[
F_\ell^{\,2}(\alpha) =\frac{q^d}{\alpha(1-\alpha)} \frac{d \alpha}{dq^d}
= \frac{1}{d}\frac{q}{\alpha(1-\alpha)} \frac{d \alpha}{dq}, \quad
F_{\ell}(\beta)^2=\frac{q}{\beta(1-\beta)} \frac{d \beta}{dq},
\]
We deduce that
\[
\frac{1}{d} \frac{d\alpha}{F_{\ell}(\alpha)^2 \alpha(1-\alpha)} =
\frac{d\beta}{F_{\ell}(\beta)^2 \beta(1-\beta)},
\]
and the relation $A_d(\alpha,\beta)=0$ is algebraic (modular equation of level $\ell$ and degree $d$). If we define the multiplier as
\[
m(\alpha, \beta)=\frac{F_{\ell}(\alpha)}{F_{\ell}(\beta)}.
\]
then
\begin{equation}\label{multiplicador}
m(\alpha, \beta)=\sqrt{\frac{1}{d} \frac{\beta(1-\beta)}{\alpha(1-\alpha)} \frac{d \alpha}{d \beta}}=\sqrt{\frac{1}{d} \frac{\beta(1-\beta)}{\alpha(1-\alpha)} \frac{\alpha'}{\beta'}}.
\end{equation} 
Hence, we have the following algebraic transformation:
\[
F_{\ell}(\alpha)=m(\alpha, \beta) F_{\ell}(\beta), \quad A_d(\alpha, \beta)=0,
\]
corresponding to level $\ell$ and degree $d$. 
From the formula for the multiplier we get
\begin{equation}\label{der-m}
\frac{m'}{\alpha'} = \frac{m}{2 \alpha'} \left( \frac{\beta'}{\beta}-\frac{\beta'}{1-\beta} -\frac{\alpha'}{\alpha} + \frac{\alpha'}{1-\alpha} + \frac{\alpha''}{\alpha'} - \frac{\beta''}{\beta'}\right).
\end{equation}
For arriving at it first square the multiplier, then take the logarithm and finally differentiate it.

\section{Explicit general formulas}

In this section we prove explicit formulas for the parameters of the Ramanujan-type series for $1/\pi$ which are more general than those derived in \cite{Gui-meth-rama} and used in \cite{Gui-rama-3A23} and \cite{Gui-Chud} because they are also valid for series with complex values of the parameters. 
\begin{theorem}
If $\beta_0 = 1-\alpha_0$, then the following formulas hold
\[
z_0=4\alpha_0 \beta_0, \qquad b=\frac{1-2\alpha_0}{\sqrt{\ell}} \left( m_0 d + \frac{1}{m_0} \right), \quad a=-2 \alpha_0 \beta_0 \frac{m'_0}{\alpha'_0} \frac{d}{\sqrt{\ell}}
\]
where
\[
m_0=\sqrt{\frac{1}{d_0} \frac{\alpha'_0}{\beta'_0}}, \quad  \tau_0=\frac{d_0 m_0}{\sqrt{\ell}} i, \quad d_0=\frac{1}{|m_0|^2}=\ell |\tau_0|^2, \quad q_0=e^{2 \pi i \tau_0}.
\]
and
\[
\frac{\beta'_0}{\alpha'_0}=\frac{1}{d m_0^2}, \qquad \frac{m'_0}{\alpha'_0} = \frac12 \left( m_0 + \frac{1}{d m_0} \right) \frac{\alpha_0-\beta_0}{\alpha_0 \beta_0} + \frac{m_0}{2\alpha'_0}\left( \frac{\alpha''_0}{\alpha'_0} - \frac{\beta''_0}{\beta'_0}\right).
\]
\end{theorem}
\begin{proof}
Apply the operator
\[
\left. a+b z_0 \frac{d}{dz} \right|_{z_0} = \left. a+ b \frac{\alpha_0 \beta_0}{1-2\alpha_0} \frac{d}{d \alpha} \right|_{u_0}
\]
to the Clausen's identity to obtain
\begin{multline}
\sum_{n=0}^{\infty} \frac{\left(\frac12\right)_n  \left(\frac1s \right)_n \left(1-\frac1s \right)_n}{(1)_n^3} z_0^n (a+bn) \\ 
= a F_{\ell}(\alpha_0) F_{\ell}(\alpha_0) 
+ \frac{b \beta_0}{1-2\alpha_0}(1+C) F_{\ell}(\alpha_0) G_{\ell}(\alpha_0) 
+  \frac{b \beta_0}{1-2\alpha_0}(1-C) F_{\ell}(\alpha_0) G_{\ell}(\alpha_0),
\end{multline}
where we intentionally have introduced the two terms involving $C$ (observe that they cancel). Then, we use the algebraic transformation
\[
F_{\ell}(\alpha)=m(\alpha, \beta) F_{\ell}(\beta),
\]
and its derivative
\[
G_{\ell}(\alpha)= \alpha  \frac{m'}{\alpha'} F_{\ell}(\beta)+\alpha \frac{m}{\beta} \frac{\beta'}{\alpha'} G_{\ell}(\beta).
\]
to substitute one factor $F_{\ell}(\alpha_0)$ of the first summand, the factor $F_{\ell}(\alpha_0)$ of the second summand, and the factor $G_{\ell}(\alpha_0)$ of the third summand. Then, identifying the coefficients of 
\[
F_{\ell}(\alpha_0) F_{\ell}(\beta_0), \quad F_{\ell}(\alpha_0) G_{\ell}(\beta_0) \quad
F_{\ell}(\beta_0) G_{\ell}(\alpha_0)
\]
to $0$, $2\alpha_0/\sqrt{\ell}$ and $2\beta_0/\sqrt{\ell}$ respectively,
we arrive at the following system of equations:
\begin{align}
a m_0 + \frac{b \beta_0 \alpha_0}{1-2\alpha_0} + \frac{m'_0}{\alpha'_0}(1-C) & = 0, \\
\frac{b m_0}{1-2\alpha_0} (1+C) & = \frac{2}{\sqrt{\ell}} \\
\frac{b m_0}{1-2\alpha_0} \frac{\beta'_0}{\alpha'_0} (1-C) & = \frac{2}{\sqrt{\ell}}.
\end{align}
From the above system we obtain the stated explicit general formulas.
\end{proof}

\section{Modular equations in Russell form}

Here, we will write the modular equations in the R. Russell's form that uses two auxiliary variables $u$ and $v$:
\[
u^h=\alpha \beta, \quad v^h=(1-\alpha)(1-\beta), \qquad P(u, v)=0,
\]
The procedure for finding $P(u,v)$ is based on the following results:

\subsection{For level 4} (Russell's Theorem) \cite[Theorem 2.1]{Chan-Liaw}. \label{teo-1}
Let $p$ be an odd prime, and write
\[
\frac{p+1}{8}=\frac{N}{D},
\]
where $N$ and $D$ are respectively the numerator and denominator of the fraction once we have simplified it. Then, if we let $h=8/D$, there exists a polynomial $P(u,v)$ of degree $N$ in $u$ and $v$ such that $P(u,v)=0$.

\subsection{For level 3} (Chan--Liaw's Theorem) \cite[Theorem 3.1]{Chan-Liaw}. \label{teo-2}
Let $p>3$ be a prime, and write
\[
\frac{p+1}{3}=\frac{N}{D},
\]
where $N$ and $D$ are respectively the numerator and denominator of the fraction after we have simplified it. Then, if we let $h=6/D$, there exists a polynomial $P(u,v)$ of degree $N$ in $u$ and $v$ such that $P(u,v)=0$.

\subsection{For level $2$}
Inspired by the above theorems, if $p$ is an odd prime number and we write
\[
\frac{p+1}{4}=\frac{N}{D},
\]
where $N$ and $D$ are respectively the numerator and denominator of the fraction after we have simplified it, and let $h=4/D$, we have observed experimentally that there exists a polynomial $P(u,v)$ of degree $(p+1)/2$ in $u$ and $v$ such that $P(u,v)=0$. \rm Anyway, when we find such a polynomial is sure that it is correct. 
  
\section{Proof of formula \ref{Ramanujan-1103} \\ using the Russell modular polynomial for $\ell=2$ and $d=29$}
We see that the formula (\ref{Ramanujan-1103}) has $z_0=99^{-4}$. As $z_0=4\alpha_0\beta_0$ where $\beta_0=1-\alpha_0$, we get
\[
\alpha_0=\frac12 - \frac{910\sqrt{29}}{9801}, \quad \beta_0=1-\alpha_0 = \frac12 + \frac{910\sqrt{29}}{9801},
\]
Then, with an approximation of $20$ digits, we have
\[
m_0=\frac{F_2(\alpha_0)}{F_2(\beta_0)} \simeq 0.1856953381770518631.
\]
and we identify $1/m_0^2 =29$. Hence, for proving (\ref{Ramanujan-1103}), we need a modular equation of degree $d=29$ for the level $\ell=2$. As we know the degree of the polynomial we get its coefficients from the mirror map by solving a linear system of equations. By doing it, we get the following modular equation:
\[
u^2=\alpha \beta, \quad v^2=(1-\alpha)(1-\beta), \qquad P(u,v)=0,
\]
where $P(u,v)$ is a very long (it is not a drawback for a computer!) symmetric polynomial of degree $15$ in $u$ and $v$ (see it in the Appendix). For applying our method we first make $\beta=1-\alpha$, which implies $u^2=v^2$. Then, if we choose $v=u$, we get a polynomial $P(u,u)$, which factors as
\[
(19602u-1)(198u+1)^2(6u+1)^2(12996u^2-3300u+1)^2(8664u^3+13164u^2+12882u+1)^2.
\]
If we choose the following root:
\[
u_0=v_0=\frac{1}{19602},
\]
then, from 
\begin{equation}\label{al-be}
u^2=\alpha \beta, \quad u^2-v^2+1 = \alpha + \beta,
\end{equation}
we get
\[
\alpha_0=\frac12-\frac{910}{9801} \sqrt{29}, \quad \beta_0=\frac12+\frac{910}{9801} \sqrt{29}.
\]
Differentiating $P(u, v)=0$ with respect to $u$ at $u=u_0$, we find
\[
v'_0=-1.
\]
Then, differentiating $P(u, v)=0$ twice with respect to $u$ at $u=u_0$, we get
\[
v''_0=\frac{352119040}{9801}.
\]
Differentiating (\ref{al-be}) with respect to $u$ at $u=u_0$, we obtain
\[
\alpha'_0=\beta'_0=\frac{1}{9801}.
\]
Then, differentiating (\ref{al-be}) twice with respect to $u$ at $u=u_0$ we obtain 
\[
\alpha''_0=-\frac{176059520}{96059601}+\frac{67432}{94743}\sqrt{29}, \quad
\beta''_0=-\frac{176059520}{96059601}-\frac{67432}{94743}\sqrt{29}
\]
and
\[
m_0=\sqrt{\frac{1}{29} \frac{\alpha'_0}{\beta'_0}}=\frac{1}{\sqrt{29}}.
\]
We obtain
\[
m'_0=-\frac{8824}{29}.
\]
Finally, using the general formulas with $\ell=2$ and $d=29$, we obtain
\[
z_0=\frac{1}{99^4}, \qquad b=\frac{52780}{9801} \sqrt{2}, \qquad a=\frac{2206}{9801} \sqrt{2},
\]
which are the values corresponding to the formula \ref{Ramanujan-1103}.

\section{Modular equations in Weber form}

Let 
\[
F_\ell(x)=\sum_{n=0}^{\infty} \frac{\left( \frac12 \right)_n \left( \frac1s \right)_n \left( 1-\frac1s \right)_n}{(1)_n^3} Z_\ell(x)^n, \qquad \ell=4 \sin^2 \frac{\pi}{s}, \qquad s=\frac{\pi}{\arcsin \frac{\sqrt{\ell}}{2}},
\]
where
\[
Z_1(x)=\frac{1728 x}{(x-16)^3}, \quad Z_2(x)=\frac{-256 x}{(x-64)^2}, \quad Z_4(x)=\frac{x}{64}.
\]
Then, the identities 
$$
F_1(x)=2 \sqrt{\frac{16-x}{64-x}} \, F_2(x)= 
\frac14 \sqrt{16-x} \, F_4(x),
$$ 
hold \cite[p. 640]{Cooper-book}, and $x=f^{24}(\tau)$, where $f(\tau)$ is the Weber function which is defined by
\[
f(\tau) = \frac{\eta^2(\tau)}{\eta \left( \frac{\tau}{2}\right)\eta(2\tau)}, \qquad 
\eta(\tau)=e^{\frac{\pi i \tau}{12}} \prod_{n=1}^{\infty} (1-e^{2\pi i \tau n}),
\]
being $\eta(\tau)$ is the Dedekind $\eta$ function.
\par For levels $\ell=1, 2, 4$ we have the following modular equations:
\begin{equation}\label{weber-modeq-d}
\alpha(1-\alpha)=Z_{\ell}(u^{24}), \quad \beta(1-\beta)=Z_{\ell}(v^{24}), \qquad \Phi_{d}(u,v)=0,
\end{equation}
where $\Phi_d(u,v)$ is the Weber polynomial of degree $d$. If $d$ is a prime $p$ these polynomials are of degree $p+1$ in $u$ and $v$, namely
\[
\Phi_d(u,v)=\sum_{i=0}^{p+1}\sum_{j=0}^{p+1} c_{i,j} u^i v^j, \quad c_{i,j}=c_{j,i}
\]
but have many coefficients that are equal to $0$ as the following theorem states:
\begin{theorem}
The coefficients $c_{i,j}$ of the Weber modular polynomials can be nonzero only if $i p  + j \equiv p+1 \pmod{24}$.   
\end{theorem}
Although these polynomials are in \cite{Weber-polynomial}, we can also determine them from the mirror map because we know the degree of them. Indeed, we only have to solve a linear system of equations which simplifies taking into account the null coefficients. Observe that the Weber polynomial of degree $29$ is much shorter than the Russell polynomial in the Appendix. In addition, we prefer to transform the polynomial $\Phi(u,v)=0$ in the following way: 
First write $Q(u^2,v^2)=uv R(u^2,v^2)$. Then, square it to get $Q^2(u^2,v^2)=u^2 v^2 R^2(u^2,v^2)$. Finally, replacing $u^2$ with $u$ and $v^2$ with $v$ we obtain $Q^2(u,v)-u v R^2(u, v)=0$. It is clear that the left hand side is a polynomial $P(u,v)$, and we have the following modular equation:
\begin{equation}\label{weber-modeq-transf-d}
\alpha(1-\alpha)=Z_{\ell}(u^{12}), \quad \beta(1-\beta)=Z_{\ell}(v^{12}), \qquad P_d(u,v)=0,
\end{equation}
We will use (\ref{weber-modeq-transf-d}) instead of (\ref{weber-modeq-d}) because the computations are simpler.

\section{Proof of formula \ref{Ramanujan-1103} \\ using the Weber modular polynomial $P_{29}(u,v)$}

\par Take $\beta_0=1-\alpha_0$. Hence, we have to find a solution of the system
\begin{equation}\label{system}
\frac{u^{12}}{(u^{12}-64)^2}=\frac{v^{12}}{(v^{12}-64)^2}, \qquad P_{29}(u,v)=0.
\end{equation}
We choose the following solution (check!):
\[
u_0=\left(\frac12-\frac{i}{2}\right) (5-\sqrt{29}), \quad 
v_0=\left(\frac12-\frac{i}{2}\right) (5+\sqrt{29})
\]
We get
\[
\alpha_0=\frac12-\frac{910}{9801} \sqrt{29}, \quad \beta_0=\frac12+\frac{910}{9801} \sqrt{29}.
\]
Differentiating $P(u, v)=0$ with respect to $u$ at $u=u_0$, we find
\[
v'_0=-\frac{27}{2}-\frac{5 \sqrt{29}}{2}.
\]
Then, differentiating $P(u, v)=0$ twice with respect to $u$ at $u=u_0$, we get
\[
v''_0=\frac{8824}{94743}(-1+i)(377+70\sqrt{29}).
\]
Differentiating (\ref{weber-modeq-transf-d}) with respect to $u$ at $u=u_0$, we obtain
\[
\alpha'_0=\beta'_0=\frac{1-i}{128079468} (5+\sqrt{29}),
\]
and we get
\[
m_0=\sqrt{\frac{1}{29} \frac{\alpha'_0}{\beta'_0}}=\frac{1}{\sqrt{29}}.
\]
Differentiating  (\ref{weber-modeq-transf-d}) twice with respect to $u$ at $u=u_0$ we obtain 
\begin{align*}
\alpha''_0 &=- \left(\frac{967391}{418435621956}+\frac{2225\sqrt{29}}{5165871876} \right) i, \\ 
\beta''_0 &=\left( \frac{1055569}{418435621956}+\frac{70339\sqrt{29}}{149810284404} \right) i.
\end{align*}
Finally, using the general formulas with $\ell=2$ and $d=29$, we obtain
\[
z_0=\frac{1}{99^4}, \qquad b=\frac{52780}{9801} \sqrt{2}, \qquad a=\frac{2206}{9801} \sqrt{2},
\]
which are the values corresponding to the formula \ref{Ramanujan-1103}.

\section{About our Maple program}

On the web page of the author there is a Maple program \cite{Gui-maple-program} that can prove automatically any Ramanujan-type series for $1/\pi$. We recommend the reader to use it to prove automatically more Ramanujan-type series for $1/\pi$, for example those in \cite{Gui-meth-rama,Gui-rama-3A23,Gui-Chud}. One of the main novelties of this paper is precisely this Maple program. For using it we introduce the level $\ell$ and the degree $d$ and the program begins computing the corresponding modular equation and continues finding the roots for $\beta=1-\alpha$ of $z_0$. In our code we have selected the roots of the polynomials that are solvable by radicals. From 
\[
z_0=4\alpha_0 (1-\alpha_0), \quad \beta_0=1-\alpha_0,
\]
we obtain
\[
\alpha_0=\frac12-\frac{1-\sqrt{1-z_0}}{2}, \qquad \beta_0=\frac12+\frac{1-\sqrt{1-z_0}}{2}.
\]
All the roots are good if we modify our formulas according to 
\[
b=\frac{1-2\alpha_0}{\sqrt{\ell}} \left( m_0 + \frac{1}{m_0 d_0} \right) d_0 \delta = 2 \frac{1-2\alpha_0}{\sqrt{\ell}} Re(m_0) d_0 \delta, \quad a=-2 \alpha_0 \beta_0 \frac{m'_0}{\alpha'_0} \frac{d_0}{\sqrt{\ell}} \delta,
\]
where 
\[
\delta=\frac{D_0}{d_0} \, \frac{{\rm Re}(M_0)}{{\rm Re} (m_0)}, 
\]
and $M_0$ the primitive multiplier and $D_0$ its corresponding degree, that is
\[
M_0=\frac{F_{\ell}(\alpha_0)}{F_{\ell}(\beta_0)}, \quad D_0=|M_0|^{-2}.
\]
On the other hand, we know that $m_0$ is the multiplier of the modular equation corresponding to the same root, and $d_0$ its degree. 
\par Finally, the program shows the corresponding Ramanujan-type series for those values of $z$ which are of the selected degree. Our code shows the times needed to obtain the modular equations and the final solutions. In our code we have tried to make the suitable simplifications but this is a difficult task and in some cases, when the simplifications are not good, the outputs are not correct. 
\par On the author's computer (2.80 GHz and 32 GB RAM) the program runs very well. For example, proving the fastest Ramanujan's formula takes $35$ seconds for obtaining automatically the modular polynomial, and $638$ seconds for proving the formula and $9$ other more.

\section*{aknowledgements}
Many thanks to Wadim Zudilin for sharing with me a generalization of the Legendre's relation. To Alin  Bostan for his great idea of simplifying using an algebraic rule. I am also very grateful  to Drew Sutherland who answered my question about which coefficients of the Weber polynomials can be nonzero. To Doron Zeilberger, who invited me give a ZOOM lecture about the Maple program cited in this paper, to Jorge Zúñiga for checking the program with many examples and detecting an error that had ocurred very few times, and to the anonymous referee for very useful information.

\newpage

\section*{Appendix I: \\ Russell Modular polynomial for level $2$ and degree $29$}

\begin{verbatim}

(u^15+v^15)
+7592191322114338383*(u^14*v+u*v^14)
-13966622568597353694807*(u^13*v^2+u^2*v^13)
+9843439764837190416529735*(u^12*v^3+u^3*v^12)
-88050856195217696119713579*(u^11*v^4+u^4*v^11)
+264881102454684464005109883*(u^10*v^5+u^5*v^10)
-339110438936155583303166131*(u^9*v^6+u^6*v^9)
+327822432596177883480907299*(u^8*v^7+u^7*v^8)
-15*(u^14+v^14)
+44151100393296850926*(u^13*v+u*v^13)
-80115946613390952012885*(u^12*v^2+u^2*v^12)
+8651254491779851942902828*(u^11*v^3+u^3*v^11)
-49687685907032270183445927*(u^10*v^4+u^4*v^10)
+7325412196773386395291122*(u^9*v^5+u^5*v^9)
+912734483936820578936054283*(u^8*v^6+u^6*v^8)
-749739327649930298909414424*(u^7*v^7)
+105*(u^13+v^13)
+111337818333508248789*(u^12*v+u*v^12)
-74674720936182298092930*(u^11*v^2+u^2*v^11)
-13575976708893635268568626*(u^10*v^3+u^3*v^10)
+175519360241020438387389123*(u^9*v^4+u^4*v^9)
-22578235036725232913012001*(u^8*v^5+u^5*v^8)
+685632655975243084716955860*(u^7*v^6+u^6*v^7)
-455*(u^12+v^12)
+159561487140115731244*(u^11*v+u*v^11)
+174695780104364520734642*(u^10*v^2+u^2*v^10)
-11101813741117032499783812*(u^9*v^3+u^3*v^9)
+201075954132759831412608055*(u^8*v^4+u^4*v^8)
+349154723555304342205053272*(u^7*v^5+u^5*v^7)
+120037255935309848803783612*(u^6*v^6)
+1365*(u^11+v^11)
+143135164756761661287*(u^10*v+u*v^10)
+283824735558969676013571*(u^9*v^2+u^2*v^9)
+2841283696852964184955593*(u^8*v^3+u^3*v^8)
-5667068077333965779571054*(u^7*v^4+u^4*v^7)
+191014637062060944039908838*(u^6*v^5+u^5*v^6)
-3003*(u^10+v^10)
+83485085314001444082*(u^9*v+u*v^9)
-14149343647392627703263*(u^8*v^2+u^2*v^8)
+1959056070055873652195160*(u^7*v^3+u^3*v^7)
-97551625043690470733963622*(u^6*v^4+u^4*v^6)
-46531375515452749172936340*(u^5*v^5)
+5005*(u^9+v^9)
+31881828039654884597*(u^8*v+u*v^8)
-167490774279805805587756*(u^7*v^2+u^2*v^7)
+2819912808713121841341252*(u^6*v^3+u^3*v^6)
-26167312382917423680651674*(u^5*v^4+u^4*v^5)
-6435*(u^8+v^8)
+7844529197113035240*(u^7*v+u*v^7)
-30988952407873752960468*(u^6*v^2+u^2*v^6)
+3116654341158724299627864*(u^5*v^3+u^3*v^5)
+1924699491602248664685678*(u^4*v^4)
+6435*(u^7+v^7)
+1196071300121674677*(u^6*v+u*v^6)
+20439274318981208511135*(u^5*v^2+u^2*v^5)
+712650893196916561319241*(u^4*v^3+u^3*v^4)
-5005*(u^6+v^6)
+105608490494850034*(u^5*v+u*v^5)
-1295519266509406265411*(u^4*v^2+u^2*v^4)
+16061879206164618424444*(u^3*v^3)
+3003*(u^5+v^5)
+4830887024610855*(u^4*v+u*v^4)
+14111349703986483150*(u^3*v^2+u^2*v^3)
-1365*(u^4+v^4)
+94634825344812*(u^3*v+u*v^3)
-17548052640024318*(u^2*v^2)
+455*(u^3+v^3)
+549304034965*(u^2*v+u*v^2)
-105*(u^2+v^2)
+368941806*(u*v)
+15*(u+v)
-1;

\end{verbatim}

\newpage
\section*{Appendix II: \\ TABLES OF RATIONAL RAMANUJAN SERIES}

In the following tables we show all the rational Ramanujan-type series for $1/\pi$, and the sequences of primitive degrees. They come from the primitive value of $\tau = {\rm Re}(\tau)+{\rm Im}(\tau) i$. Observe that adding an integer $n$ then the value of $q=e^{2 \pi i \tau}$ is preserved. Hence we get the sequence of degrees: $d_n=\ell ((n+{\rm Re}(\tau))^2+{\rm Im}(\tau)^2)$. There exists also other possible sequences of non-primitive degrees corresponding to some $\delta \neq 1$. We can discover these other sequences using our program (remember that it is only prepared for working with prime numbers).

\begin{table}[ht]
    \begin{tabular}{|c c c c | c c c c|}
        \hline &&&&&&& \\ [-1ex]
        \rm Primitive $d_n$ & $a$ & $b$ & $z<0$ &\rm Primitive $d_n$ & $a$ & $b$ & $z>0$ \\ &&&&&&& \\ [-1ex]
        \hline \hline &&&&&&& \\
        $n^2+n+2$ & $\frac{8}{5\sqrt{15}}$ & $\frac{63}{5\sqrt{15}}$ & $-\frac{4^3}{5^3}$ &
        $n^2+2$ & $\frac{3}{5\sqrt5}$ & $\frac{28}{5\sqrt5}$ & $\frac{3^3}{5^3}$ \\ &&&&&&& \\
        $n^2+n+3$ & $\frac{15}{32\sqrt2}$ & $\frac{154}{32\sqrt2}$ & $-\frac{3^3}{8^3}$ &
        $n^2+3$ & $\frac{6}{5\sqrt{15}}$ & $\frac{66}{5\sqrt{15}}$ & $\frac{4}{5^3}$ \\ &&&&&&& \\
        $n^2+n+5$ & $\frac{25}{32\sqrt6}$ & $\frac{342}{32\sqrt6}$ & $-\frac{1}{8^3}$ &
        $n^2+4$ & $\frac{20}{11\sqrt{33}}$ & $\frac{252}{11\sqrt{33}}$ & $\frac{2^3}{11^3}$ \\ &&&&&&& \\
        $n^2+n+7$ & $\frac{279}{160\sqrt{30}}$ & $\frac{4554}{160\sqrt{30}}$ & $-\frac{9}{40^3}$ &
        $n^2+7$ & $\frac{144\sqrt3}{85\sqrt{85}}$ & $\frac{2394\sqrt{3}}{85\sqrt{85}}$ & $\frac{4^3}{85^3}$ \\ &&&&&&& \\
        $n^2+n+11$ & $\frac{526\sqrt{15}}{80^2}$ & $\frac{10836 \sqrt{15}}{80^2}$ & $-\frac{1}{80^3}$ &&&& \\ &&&&&&& \\
        $n^2+n+17$ & $\frac{10177\sqrt{330}}{3 \cdot 440^2}$ & $\frac{261702\sqrt{330}}{3\cdot 440^2}$ & $-\frac{1}{440^3}$ &&&& \\ &&&&&&& \\
        $n^2+n+41$ & $\frac{27182818\sqrt{10005}}{3 \cdot 53360^2}$ & $\frac{1090280268\sqrt{10005}}{3\cdot 53360^2}$ & $-\frac{1}{53360^3}$ &&&& \\ &&&&&&& \\
        \hline
    \end{tabular}
    \vskip 0.25cm
    \caption{Rational Ramanujan-type series of $\ell=1$ for $1/\pi$}
\end{table}

\begin{table}[b]
    \begin{tabular}{|c c c c | c c c c|}
        \hline &&&&&&& \\ [-1ex]
        \rm Primitive $d_n$ & $a$ & $b$ & $z<0$ & \rm Primitive $d_n$ & $a$ & $b$ & $z>0$ \\ &&&&&&& \\ [-1ex]
        \hline \hline &&&&&&& \\ 
        $2n^2+n+1$ & $\frac{1}{\sqrt3}$ & $\frac{5}{\sqrt3}$ & $-\frac{16}{9}$ &
        $2n^2+1$  & $\frac{8}{18} i$ & $\frac{35}{18} i$ & $\frac{256}{81}$ \\ &&&&&&& \\
        $2n^2+2n+3$ & $\frac{3}{8}$ & $\frac{20}{8}$ & $-\frac{1}{4}$ &
        $2n^2+2$ & $\frac{2}{9}$ & $\frac{14}{9}$ & $\frac{32}{81}$ \\ &&&&&&& \\
        $2n^2+2n+4$ & $\frac{8}{9\sqrt7}$ & $\frac{65}{9\sqrt7}$ & $-\frac{16^2}{63^2}$ &
        $2n^2+3$ & $\frac{1}{2\sqrt3}$ & $\frac{8}{2\sqrt3}$ & $\frac{1}{9}$ \\ &&&&&&& \\
        $2n^2+2n+5$ & $\frac{3\sqrt3}{16}$ & $\frac{28\sqrt3}{16}$ & $-\frac{1}{48}$ &
        $2n^2+5$ & $\frac{4}{9\sqrt2}$ & $\frac{40}{9\sqrt2}$ & $\frac{1}{81}$ \\ &&&&&&& \\
        $2n^2+2n+7$ & $\frac{23}{72}$ & $\frac{260}{72}$ & $-\frac{1}{18^2}$ &
        $2n^2+9$ & $\frac{27}{49\sqrt3}$ & $\frac{360}{49\sqrt3}$ & $\frac{1}{7^4}$ \\ &&&&&&& \\
        $2n^2+2n+13$ & $\frac{41\sqrt5}{288}$ & $\frac{644\sqrt5}{288}$ & $-\frac{1}{5 \cdot 72^2}$ &
        $2n^2+11$ & $\frac{19}{18\sqrt{11}}$ & $\frac{280}{18\sqrt{11}}$ & $\frac{1}{99^2}$ \\ &&&&&&& \\
        $2n^2+2n+19$ & $\frac{1123}{3528}$ & $\frac{21460}{3528}$ & $-\frac{1}{882^2}$ &
        $2n^2+29$ & $\frac{4412}{9801\sqrt2}$ & $\frac{105560}{9801\sqrt2}$ & $\frac{1}{99^4}$ \\ &&&&&&& \\
        \hline
    \end{tabular}
    \vskip 0.25cm
    \caption{Rational Ramanujan-type series of $\ell=2$ for $1/\pi$}
\end{table}

\begin{table}[ht]
    \begin{tabular}{|c c c c | c c c c|}
        \hline &&&&&&& \\ [-1ex]
        \rm Primitive $d_n$ & $a$ & $b$ & $z<0$ & \rm Primitive $d_n$ & $a$ & $b$ & $z>0$ \\ &&&&&&& \\ [-1ex]
        \hline \hline &&&&&&& \\
        $3n^2+3n+1$ & $\frac{4}{3\sqrt3}$ & $\frac{15}{3\sqrt3}$ & $-4$ &
        $-$ & $\frac{3}{10} i$ & $\frac{10}{10} i$ & $\frac{27}{2}$ \\ &&&&&&& \\
        $3n^2+3n+3$ & $\frac{\sqrt{3}}{4}$ & $\frac{5\sqrt{3}}{4}$ & $-\frac{9}{16}$ &
        $-$ & $\frac{3}{12} i$ & $\frac{11}{12} i$ & $\frac{27}{16}$ \\ &&&&&&& \\
        $3n^2+3n+5$ & $\frac{7}{12\sqrt3}$ & $\frac{51}{12\sqrt3}$ & $-\frac{1}{16}$ &
        $3n^2+2$ & $\frac{1}{3\sqrt{3}}$ & $\frac{6}{3\sqrt{3}}$ & $\frac{1}{2}$ \\ &&&&&&& \\
        $3n^2+3n+7$ & $\frac{\sqrt{15}}{12}$ & $\frac{9\sqrt{15}}{12}$ & $-\frac{1}{80}$ &
        $3n^2+4$ & $\frac{8}{27}$ & $\frac{60}{27}$ & $\frac{2}{27}$ \\ &&&&&&& \\
        $3n^2+3n+11$ & $\frac{106}{192\sqrt3}$ & $\frac{1230}{192\sqrt3}$ & $-\frac{1}{2^{10}}$ &
        $3n^2+5$ & $\frac{8}{15\sqrt3}$ & $\frac{66}{15\sqrt3}$ & $\frac{4}{125}$ \\ &&&&&&& \\
        $3n^2+3+13$ & $\frac{26\sqrt7}{216}$ & $\frac{330\sqrt7}{216}$ & $-\frac{1}{3024}$
        &&&& \\ &&&&&&& \\
        $3n^2+3n+23$ & $\frac{827}{1500\sqrt3}$ & $\frac{14151}{1500\sqrt3}$ & $-\frac{1}{500^2}$
        &&&& \\ &&&&&&& \\
        \hline
    \end{tabular}
    \vskip 0.25cm
    \caption{Rational Ramanujan-type series of $\ell=3$ for $1/\pi$}
\end{table}

\begin{table}[p]
    \begin{tabular}{|c c c c | c c c c|}
        \hline &&&&&&& \\ [-1ex]
        \rm Primitive $d_n$ & $a$ & $b$ & $z<0$ & \rm Primitive $d_n$ & $a$ & $b$ & $z>0$ \\ &&&&&&& \\ [-1ex]
        \hline \hline &&&&&&& \\ 
        $4n^2+4n+1$ & $1$ & $3$ & $-8$ &
        $4n^2+3n+1$ & $\frac{8}{4} i$ & $\frac{21}{4} i$ & $64$ \\ &&&&&&& \\
        $4n^2+4n+3$ & $\frac14$ & $\frac44$ & $-1$ &
        $4n^2+2n+1$ & $\frac{1}{2} i$ & $\frac{3}{2} i$ & $4$ \\ &&&&&&& \\
        $4n^2+4n+5$ & $\frac{1}{2\sqrt2}$ & $\frac{6}{2\sqrt2}$ & $-\frac{1}{8}$ &
        $4n^2+3$ & $\frac14$ & $\frac64$ & $\frac14$ \\ &&&&&&& \\
        $$ & $$ & $$ & $$ &
        $4n^2+7$ & $\frac{5}{16}$ & $\frac{42}{16}$ & $\frac{1}{64}$ \\ &&&&&&& \\
        \hline
    \end{tabular}
    \vskip 0.25cm
    \caption{Rational Ramanujan-type series of $\ell=4$ for $1/\pi$}
\end{table}

\end{document}